\nonstopmode
\documentclass[12pt,a4paper]{article}

\usepackage{amsmath,amsthm, amssymb}
\usepackage{latexsym}
\usepackage{epsfig}
\usepackage{url,paralist}
\usepackage{bbm}
\newcommand{\cc}{\mathcal{C}}
\newcommand{\stab}{\mbox{\rm STAB}(G)}
\newcommand{\gstab}{G_{\mbox{\rm \scriptsize STAB}}}
\newcommand{\reals}{\mathbbm{R}}
\newcommand{\suchthat}{\, | \,}

\newtheorem{lemma}{Lemma}[section]
\newtheorem{theorem}[lemma]{Theorem} 
\newtheorem{defin}[lemma]{Definition}

\begin{document}

\title{Edge Expansion of Cubical Complexes}
\author{Thomas Voigt\thanks{%
Supported by Deutsche Forschungs-Gemeinschaft (DFG),
Forschergruppe ``Algorithms, Structure, Randomness''}\\
\small Inst.\ Mathematics, MA 6-2\\[-1.2mm]
\small TU Berlin, D-10623 Berlin, Germany\\[-1.2mm]
\small\url{tvoigt@math.tu-berlin.de}}
\date{}

\maketitle

\begin{abstract}
In this paper we show that graphs of ``neighbourly'' cubical complexes 
-- cubical complexes in which every pair of vertices
spans a (unique) cube -- have good expansion properties, using a technique based on 
multicommodity flows. By showing that graphs of stable
set polytopes are graphs of neighbourly cubical complexes we give a new proof
that graphs of stable set polytopes have edge expansion 1.

\end{abstract} 

\section{Introduction}
$0/1$-polytopes arise naturally in a great variety of 
contexts, yet only few combinatorial properties 
are well understood (for an introduction see \cite{Zie00}). 
Special attention was paid to polytope classes that correspond to
some key combinatorial optimization problems.  \\

One of these properties is the edge expansion of
the graph of $0/1$-polytopes, where the edge expansion $\chi(G)$ 
of a graph $G=(V,E)$ is defined as
\[ \chi(G) = \min_{ X \subset V \atop 0 < |X| \leq \frac{1}{2}|V|} \frac{|\delta(X)|}{|X|} \]
where $\delta(X)$ is the set of edges with one endpoint in $X$ and the other
endpoint in $V \setminus X$.
Mihail and Vazirani conjectured $\chi(G) \geq 1$ for all graphs 
arising as $1$-skeleton of $0/1$-polytopes \cite{Kaibel03}\cite{Mih92}. This
lower bound was verified for a number of polytope classes, among them
matroid basis polytopes of balanced matroids $\cite{FM92}$ and graphs of stable set
polytopes \cite{Kaibel03}\cite{Mih92}. Lower bounds on the expansion of the 1-skeleton 
of polytopes are of particular interest because they correspond to
upper bounds on the mixing time of the natural random
walk on the graph of the polytope (where a transition between
two vertices has positive probability if and only if the
vertices are connected by an edge). This was successfully
used to obtain algorithms for random sampling and approximate
counting of the vertices in the case of balanced matroids \cite{JS02}. \\
  
In this paper we will consider the edge expansion of 
neighbourly (abstract) cubical complexes:
\begin{defin}
An {\bf \emph{(abstract) cubical complex}} is a pair $K = (V,\cc)$ of a set $V$ and 
a subset $\cc$ of the power set of $V$ (``cubes'') such that 
\begin{enumerate}
\item $c_1, c_2 \in \cc \Rightarrow c_1 \cap c_2 \in \cc$
\item For every $c \in \cc$ there is a $d \geq 0$ and a bijection $\phi_c$ between
$\{0,1\}^d$ and $c$ such that the image $\phi_c(s)$ of a subset $s \subset \{0,1\}^d$ is 
is a member of $\cc$ if and only if $s$ is a face of $\{0,1\}^d$.
\end{enumerate}
A graph $G=(V,E)$ is the graph of an abstract cubical complex if
there is a cubical complex $(V,\cc)$ such that $\{v_1, v_2\} \in E \Leftrightarrow \{v_1, v_2\} \in \cc$. \\

A cubical complex is {\bf \emph{neighbourly}} if for every pair of vertices $v, w$ 
there is a  cube $c \in \cc$
such that $v$ and $w$ are vertices in $c$. 
\end{defin}
If $G$ is the graph of a cubical complex then the induced subgraph $G[c]$ 
is a cube graph for every $c \in \cc$. Cubes have edge expansion 1,
which extends to neighbourly cubical complexes: 
\begin{theorem}
\label{thm:complexexpansion}
Let $G=(V,E)$ be the graph of a neighbourly cubical complex $K=(V,\cc)$. 
Then the edge expansion of $G$ satisfies $\chi(G) \geq 1$.
\end{theorem}
We will prove Theorem \ref{thm:complexexpansion} in Section
\ref{sec:ccubes}. \\
In Section \ref{sec:stable} we will turn to the structure of  
\emph{stable set polytopes}, where 
the stable set polytope $\stab$ of a graph $G = (V,E)$ is
defined as the convex hull of the incidence vectors of
all stable sets in $G$. $\stab$ Every stable
set polytope $\stab$ corresponds to a neighbourly cubical
complex $(V(\gstab), \cc)$ in a very natural way, which yields 
a new proof for the edge expansion of stable set polytopes. 
($\gstab$ is the graph of $\stab$.)\\

\begin{figure}[ht]
\label{fig:gstab}
\begin{center}
\epsfig{file=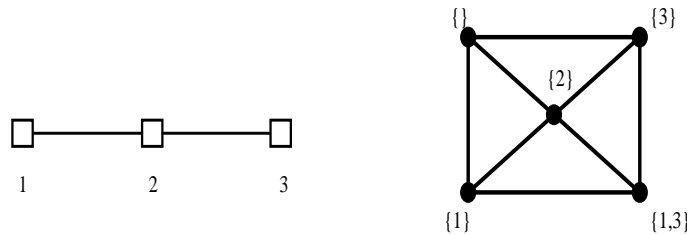,width=9cm, height=3cm}
\caption{$P_3$ and the graph of $\gstab(P_3)$}
\end{center}
\end{figure}

However, the natural random walk on stable set polytopes cannot
be used to construct an approximate counting algorithm for the
number of stable sets of a graph $G$: Even though there
exist fully polynomial approximation algorithms if $G$ has maximal degree up to 4  
using a different random walk \cite{LV97}\cite{DG99}, Dyer, Frieze and 
Jerrum \cite{DFJ02} proved that it is not possible 
to approximate the number of stable sets in polynomial
time if the degree of $G$ is larger than~25, unless RP=NP. Since
$\chi(\gstab) \geq 1$ for any graph $G$ we know
that a random walk on this
graph would be mixing in time polynomial in $\log(|V(\gstab)|)$ and
it follows from the previous remarks that it is not possible to 
use this neighbourhood structure of $\gstab$ for a 
rapidly mixing random walk -- 
which isn't surprising since the number of neighbours of a 
stable set (in this neighbourhood structure) can be exponential in 
the size of $G$, which makes it difficult to design a random walk $X$ that
 converges to the uniform distribution in polynomial time.

\section{Edge Expansion of cubical complexes}
\label{sec:ccubes}

In order to find a lower bound on the edge expansion of neighbourly 
cubical complexes we use a flow
method that may be viewed as a generalization of canonical paths
(see \cite{JS96}) and was made explicit in \cite{Kaibel03} \cite{MS04}:
Define for each ordered pair of vertices $v, w \in V$
a flow $f_{vw}: E \rightarrow \reals$ that sends one unit of flow
from $v$ to $w$. We 
obtain the desired result by giving an upper bound on the total directed flow 
$F := \sum_{v,w \in V} f_{vw}$ on every edge: Suppose the maximal
directed flow on any edge is $\mu$. For every subset $X \subset V, |X| \leq \frac{1}{2} |V|$ we have by
definition a total flow from $X$ to $V\setminus X$ of at least $|X| |V \setminus X|$. Since
every cut edge contributes at most $\mu$ to this flow we have $ |\delta(X)| \mu \geq |X||V \setminus X| $
and therefore 
\[ \frac{|\delta(X)|}{|X|} \geq \frac{|V \setminus X|}{\mu} \geq \frac{|V|}{2 \mu}. \]
In particular, if $\mu \leq \frac{1}{2}|V|$ then $\chi(G) \geq 1$. \\

It will be convenient to consider $G$ as a 
directed graph (every undirected edge is replaced by two
directed edges), so we just have to give a bound for the
maximal flow on any directed edge.

Since $G$ is the graph of a cubical complex there is a canonical
way to define a flow $f_{vw}$ between two vertices $v,w$ that are in some cube
$c \in \cc$:  Let $c_{vw} \in \cc$ be the
subcube spanned by $v$ and $w$ (so $v$ and $w$ are antipodal in $c_{vw}$) and distribute the
flow equally over all shortest paths from $v$ to $w$ in the (directed) cube graph $G[c_{vw}]$. \\

Now Theorem \ref{thm:complexexpansion} follows easily from a lemma that 
may be interesting in its own right:
\begin{lemma} 
\label{lem:complexlemma}
Let $K=(V, \cc)$ be any abstract cubical complex 
and $e$ a $1$-dimensional cube of the complex. Let 
$C_e = \{c \in \cc \suchthat e \subset c\}$ be the set of cubes that contain $e$.
Define a flow $f_{vw}$ as described above for every pair of vertices $v,w \in V$ 
that span a cube in $\cc$ and let $F=\sum_{v,w \in V} f_{vw}$ (where $f_{vw}$ is
zero if there is no cube $c$ such that $v,w \in c$).  Then
\[ F(e) = |C_e|  \leq \tfrac{1}{2} |V|.\]
\end{lemma} 
{\bf Proof:}  
First observe that for any $d$-dimensional cube $c \in \cc$ the flow (on directed edges)
\[ f_c(e) := \sum_{v,w \in c} f_{vw}(e) \] 
induced by antipodal pairs of vertices in $c$ is sent using edges of $G[c]$ and
adds up to $1$ for every edge $e \in G[c]$. 
To see this, check that there are $2^{d-1}$ antipodal pairs of vertices
and every pair generates a total flow of $2d$, so the flow on each 
of the $d2^d$ edges is 1. So we have
\[ F(e) = \sum_{c \in C_e} f(c) = |C_e|. \] 
But the number of cubes in $C_e$ is limited by $\frac{1}{2}|V|$: 
Given any $c \in C_e$ there exists an unique edge $e_a(c)$ that is
antipodal to $e$ in $c$. Since $e$ and $e_a(c)$ span $c$
($c$ is the unique inclusion-minimal cube that contains both edges)
we know that different cubes $c, \tilde c$ containing $e$ are mapped to
vertex-disjoint different edges $e_a(c), e_a({\tilde c})$. 
So the lemma follows from the fact that there exist $|C_e|$ 
pairwise vertex-disjoint edges in $\cc$. \\

{\bf Proof of Theorem \ref{thm:complexexpansion}:} 
Since every pair of vertices spans a cube, the flow constructed in Lemma
\ref{lem:complexlemma} satisfies that every vertex sends one commodity of flow to
every other vertex. The flow on any edge $e$ is at most
$\frac{1}{2} |V|$ and therefore the edge expansion of $G$ is 
at least 1. \begin{flushright} $\Box$ \end{flushright}

\section{Graphs of stable set polytopes}
\label{sec:stable}

The idea of using cubic subgraphs for expansion of 0/1 polytopes was used
by Milena Mihail in \cite{Mih92} to show that the perfect matching polytope
of any graph has edge expansion 1. Volker Kaibel generalized this concept in \cite{Kaibel03}
to polytopes with ``cube-spanned walls'' which includes stable set polytopes and matching polytopes.  
The proof presented here is similar but exploits the fact that
the graph of a stable set polytope is the graph of a neighbourly cubical complex. \\

\begin{lemma}
\label{lemmacc}
For any graph $G$ there exists a neighbourly cubical
complex $K=(V(\gstab), \cc)$ such that $\gstab$ is the graph of
$K$.
\end{lemma}

Every arbitrary graph $H=(V_H, E_H)$ corresponds to the trivial  
cubical complex $K_H = (V_H, \cc)$ with $\cc = V_H \cup \{ \{x, y\} | (x,y) \in E_H\}$. 
This complex may be extended to a neighbourly cubical complex if~$H$
can be covered by a set of suitable cube graphs that is closed
under taking intersections.
\begin{figure}[ht]
\label{fig:ccomplex}
\begin{center}
\epsfig{file=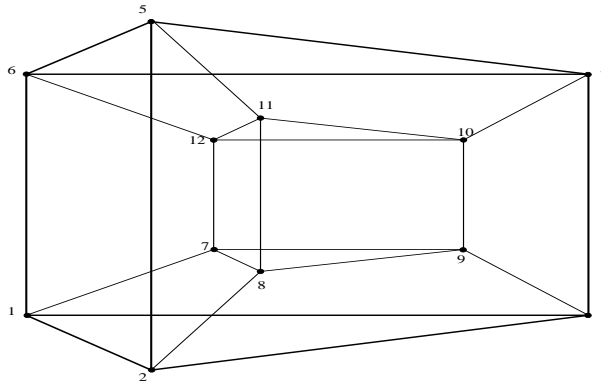,width=8cm, height=5cm}
\caption{A graph of a neighbourly cubical complex}
\end{center}
\end{figure}
Figure \ref{fig:ccomplex} shows  the graph
of the $4$-dimensional polytope $P$ arising as the product of a triangle and a 
square quadrilateral. This graph on 12 vertices that can be decomposed as the union
of the 3-dimensional cubes  induced by the vertices $\{1,2,5,6,7,8,11,12\}$,
$\{2,3,4,5,8,9,10,11\}$ and $\{1,3,4,6,7,9,10,11\}$. By adding these three cubes and all
their faces to the trivial cubical complex we obtain a neighbourly cubical complex. 
Note that $P$ is the stable set polytope of the simple graph on 4 vertices with
exactly one edge.   \\

For arbitrary stable set polytopes we can explicitly construct the cubes
for all pairs of vertices: \\

{\bf Proof of Lemma \ref{lemmacc}:} 
For 
$\tilde V \subset V$ let $G[\tilde V]$ be the induced subgraph on the
vertices $\tilde V$. 

We will define a neighbourly cubical complex such that $\gstab$ is the graph of this complex. \\

Chv\'atal  \cite{Chvatal75} proved that vertices corresponding to stable sets $s, t$ are adjacent
in $\gstab$ if and only if $G[s \vartriangle t]$ is connected. (We need only the ``if'' part
in this proof.)
Fix arbitrary stable sets $s,t$ and let $A_1, \ldots, A_k \subseteq V$ 
be the vertex sets of the connected components of the bipartite graph $G[s \vartriangle t]$.
In this graph we can choose $2^k$ maximal independent sets by choosing independently
for each $j$ either all the vertices of $A_j \cap s$ or all the 
vertices of $A_j \cap t$. Any pair of these independent sets
is connected by an edge of $\gstab$ if and only if they differ in exactly
one component, so the subgraph induced by these $2^k$ sets is the 
graph of a $k$-dimensional cube and the faces of the cube correspond
to the stable sets where we fix the choice in a number of components. \\

\begin{figure}[ht]
\label{fig:cubes}
\begin{center}
\epsfig{file=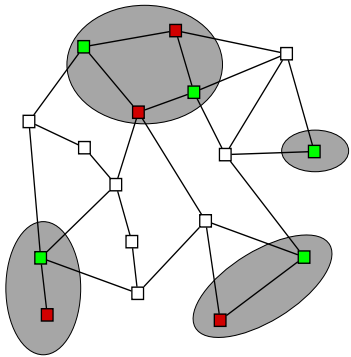,width=5cm, height=5cm}
\caption{The bipartite subgraph defined by two stable sets}
\end{center}
\end{figure}

More formally, we take the set of these $2^k$ stable sets $c_{st}$ 
as a cube in our cubical complex by defining 
\[ c_{st} := \big\{  s \vartriangle \big( \bigcup_{j \in J} A_j \big)
\suchthat J \subseteq [k] \big\} \]
where $[k] = \{1, \ldots, k\}$.
By construction any pair of stable
sets defines a cube, so we can define the neighbourly cubical complex
$K=(V(\gstab), \cc)$ where $\cc = \{ c_{st} | s,t \mbox{ stable sets of } G \}$, provided 
that the intersection of two cubes $c_1, c_2 \in \cc$ is again a cube in $\cc$. But this is the case: \\
Let $s,t$  be as above, set $\mathcal{A} := \{ A_1 \ldots A_k\}$ and pick 
another arbitrary pair of stable sets $\tilde s, \tilde t$ with sets
$\mathcal{\tilde A} := \{\tilde A_1, \ldots, \tilde A_{\tilde k}\}$ such that $G[\tilde A_j]$ are the connected
components of $G[\tilde s \vartriangle \tilde t]$ for $j \leq \tilde k$. Assume that
the intersection of the cubes $c_{st}$ and $c_{\tilde s \tilde t}$ contains at least
one stable set $v$ (otherwise there is nothing to show). Let $l$ be the number
of sets $A_j \in \mathcal{A}$ that are identical to some  $\tilde A_{\tilde j} \in \mathcal{\tilde A}$ (possibly $l=0$)
and assume w.l.o.g. $A_j = \tilde A_j$ for all $j \leq l$. We claim that
\[ c_{st} \cap c_{\tilde s \tilde t} = \big\{ v \vartriangle  \big( \bigcup_{i \in I} A_i \big) \suchthat I \subset [l] \big\} \]
which is a cube in $\cc$ since it is spanned by the stable sets $v$ and $v \vartriangle (\cup_{i \leq l} A_i)$.
Clearly $v \vartriangle  ( \bigcup_{i \in I} A_i) \in c_{st} \cap c_{\tilde s \tilde t}$ for every 
$I \subset [l]$, so we have to show that for any $w \in c_{st} \cap c_{\tilde s \tilde t}$ there exists
an index set $I \subset [l]$ such that $w = v \vartriangle  ( \bigcup_{i \in I} A_i)$.
By definition we have sets $J \subset [k], \tilde J \subset [\tilde k]$ such that
\[ w = v \vartriangle  \big( \bigcup_{j \in J} A_j \big) = v \vartriangle \big( \bigcup_{j \in \tilde J} \tilde A_j \big) \]
But since the $A_j$ (and $\tilde A_j$ respectively) induce disconnected regions of $G$,
the condition $\bigcup_{j \in J} A_j = \bigcup_{j \in \tilde J} \tilde A_j$ implies
that  there is a bijection $\phi: J \rightarrow \tilde J$ such
that $A_j = \tilde A_{\phi(j)}$ for all $j \in J$, therefore $J \subset [l]$.  \\
So $\gstab$ is the graph of the complex $V, \cc$ and by Theorem \ref{thm:complexexpansion} 
\[ \chi(\gstab) \geq 1. \]
\begin{flushright} $\Box$ \end{flushright}

{\bf Remark:} 
Let $s, t, k$ and $A_j$ for $j\leq k$ be defined as above.
Then the vertices of $c_{st}$ do span a face $F_{st}$ of
$\stab$ (see \cite{Kaibel03}, the face is defined by the 
equations $x_j=0$ for $j \in V \setminus (s \cup t)$
and $x_i = 1$ for $i \in s \cap t$). However, by choosing all
vertices in $s \cup t$ and
adding suitable subsets of $A_j \cap s$ or $A_j \cap t$ for each $j \leq k$ we can
construct stable sets that are contained in $F_{st}$ but not in $c_{st}$, so
the complex $(V(\gstab), \cc)$ is not a subcomplex of $\stab$.

\section*{Acknowledgements}
I thank Volker Kaibel for many helpful suggestions which led to
a significantly simpler proof of Lemma \ref{lem:complexlemma}.

\bibliography{/homes/combi/tvoigt/works/Bibtex/BasicBib}

\end{document}